\documentclass{article}
\usepackage{amsmath}
\usepackage{amsfonts}
\usepackage{amssymb}

\def\1{{\bf 1}}
\def\2{{\bf 2}}
\def\3{{\bf 3}}
\def\j{{\bf j}}

\def\b{\beta}

\def\D{\Delta}

\def\R{\mathbb{R}}

\def\n{\noindent}

\def\p{\partial}

\def\t#1{{\widetilde{#1}}}

\def\X{{\bf X}}
\def\Y{{\bf Y}}
\def\b#1{{\overline{#1}}}

\def\h#1{{\widehat{#1}}}
\def\u#1{{\underline{#1}}}
\def\0{\mid_{\,0}}

\newtheorem{thm}{Theorem}

\newcommand{\be}{\begin{equation}}
\newcommand{\ee}{\end{equation}}
\newcommand{\ba}{\begin{array}}
\newcommand{\ea}{\end{array}}

\title{Local moduli in the special 2-flags of length 5}

\author{Piotr Mormul\\
          \\
\small Institute of Mathematics, University of Warsaw, Poland\\
\small e-mail: mormul@mimuw.edu.pl\normalsize}
\date{}
\begin{document}
\maketitle
\begin{abstract}
A number of key issues concerning distributions generating 1-flags 
(most often called Goursat flags) has been settled over the past 30 years.\\
Presently similar questions are being discussed as regards distributions generating 
{\it multi}-flags. (More precisely, only so-called {\it special\,} multi-flags, 
to avoid functional moduli in local classifications.) In particular, special 2-flags 
of small lengths are a natural ground for the search of generalizations of theorems 
established earlier for Goursat structures.  This includes the search for the first 
appearing modulus (or moduli) in the classification up to local diffeomorphisms of 
special 2-flags. (For Goursat flags the first modulus of the local classification 
appears in length 8.) 
It has been known in this respect that up to length 4 that classification is finite 
(\cite{MP1}, \cite{CH}), and that in length 7 at least one numerical modulus exists 
(as produced on pp.\,39\,--\,41 in \cite{MP1}). In the last fully classified length 4 
possible are precisely 34 local geometries (local models) of special 2-flags.\\
We now demonstrate that in the length 5 single numerical moduli show up in 
exactly three out of altogether 41 singularity classes existing in that 
length. 
\end{abstract}
\section{Definition of special $m$-flags (encoded in Sandwich Diagram)}
Special $m$-flags (the natural parameter $m \ge 2$ is sometimes called 
the `width'; an alternative name in \cite{CH} is `Goursat multi-flags') 
of lengths $r \ge 1$ are defined as follows.\\

For a distribution $D$ on a manifold $M$, the tower of consecutive 
Lie squares of $D$ (the so-called {\it derived flag} of $D$) 
\[
D = D^r \subset D^{r-1} \subset D^{r-2} \subset \,\cdots \,\subset 
D^1 \subset D^0 = TM
\]
(that is, $[D^j,\,D^j] = D^{j-1}$ for $j = r,\,r-1,\dots,\,2,\,1$) 
should consist of distributions of ranks, starting from the smallest 
object $D^r$: $m+1$, $2m+1,\,\dots$, $rm + 1$, $(r+1)m + 1 = 
{\rm dim}\,M$ \,such \,that 
\vskip.7mm
\n$\bullet$ for $j = 1,\dots,\,r-1$ the Cauchy-characteristic module 
$L(D^j)$ of $D^j$ sits already in the smaller object $D^{j+1}$: 
$L(D^j) \subset D^{j+1}$ {\bf and} is regular of corank 1 in 
$D^{j+1}$, while $L(D^r) = 0$\,; 
\vskip1mm
\n$\bullet\bullet$ the covariant subdistribution $F$ of \,$D^1$ (see 
\cite{KR}, p.\,5 for the definition extending the classical Cartan's 
approach from \cite{C}, p.\,121) exists and is involutive. \,Note that, 
in view of Lemma 1 in \cite{KR}, such an $F$ is {\bf automatically} 
of corank 1 in $D^1$; the hypotheses in that lemma are satisfied as 
\,${\rm rk}\,[D^1,\,D^1]/D^1 \,= \,m > 1$.\footnote{\,\,Equivalently, 
using Tanaka's and Yamaguchi's approach \cite{T,Y} (well anterior to 
\cite{KR} and not designed for special flags, although applicable 
for them), one stipulates in $\bullet\bullet$ two things: ---\,the 
distribution $\h{D}^1 = D^1/L(D^1)$ of rank $m+1$ on a manifold of 
dimension $2m+1$ is of type $\frak{C}^1(1,\,m)$ of \cite{T} and, 
as such, possesses its {\it symbol\,} subdistribution 
$\h{F} \subset \h{D}^1$ (\cite{Y}, p.\,30) \,and ---\,$\h F$ is 
involutive (cf. Prop.\,1.5 in \cite{Y}). $F$ is then the counterimage 
of $\h{F}$ under the factoring out by $L(D^1)$. Thus, for special 
$m$-flags, $m \ge 2$, the stipulated involutive corank one 
subdistribution of $D^1$ is at the same time: the covariant 
subdistribution in the Cartan-Kumpera-Rubin sense {\it and\,} 
symbol subdistribution in the Tanaka-Yamaguchi sense.}\\

\n{\it Attention.} The reference works \cite{A} and \cite{SY} have 
remarkably compactified the above original definition. In the light of 
those contributions, the extensive condition $\bullet$ is altogether 
{\it redundant}. For it follows from the condition $\bullet\bullet$ 
and the property of regular dimension growth in the derived flag 
of the initial distribution $D$.\\

However, the redundancy in the original definition has its advantages, 
too. Namely, the objects entering that definition can be assembled 
into so-called {\it Sandwich Diagram}.\footnote{\,\,so called after 
a similar (if not identical) diagram assembled for Goursat distributions, 
or 1-flags, in \cite{MZ1}} 
\[
\!\!\!\!\!\!\!\!\!\!\!\!\!\!\!\ba{ccccccccccccc}
TM = D^0 & \supset & D^1 & \supset & D^2 & \supset & \cdots & \supset
& D^{r-1} & \supset & D^r & & \\
& & \cup & & \cup & & & & \cup & & \cup & & \\
& & F & \supset & L(D^1) & \supset & \cdots & \supset & L(D^{r-2})
& \supset & L(D^{r-1}) & \supset & L(D^r) = 0
\ea
\]
(Thus the version of the definition recalled in \cite{CH} on p.\,424, 
putting into relief the Sandwich Diagram itself, is misleading. 
Reiterating, the Sandwich Diagram is {\it not\,} a part of 
the definition, but a consequence of it.)\\

\n Note that the inclusions $L(D^{j-1}) \supset L(D^j)$ in the diagram's lower 
line are due to the Jacobi identity. As for the inclusion $F \supset L(D^1)$, 
it follows from \cite{KR}, too. All vertical inclusions in this diagram are 
of codimension one, while all (drawn, not meaning superpositions of them) 
horizontal inclusions are of codimension $m$. The squares built by these 
inclusions can, indeed, be perceived as certain `sandwiches'. For instance, 
in the utmost left sandwich $F$ and $D^2$ are as if fillings, while $D^1$ 
and $L(D^1)$ constitute the covers (of dimensions differing by $m+1$, 
one has to admit). At that, the sum $m+1$ of {\bf co}dimensions, in $D^1$, 
of $F$ and $D^2$ equals the dimension of the quotient space $D^1/L(D^1)$, 
so that it is natural to ask how the $m$-dimensional space $F/L(D^1)$ 
and the line $D^2/L(D^1)$ are mutually positioned in $D^1/L(D^1)$. 
Similar questions impose by themselves in further sandwiches `indexed' 
by the upper right vertices $D^3,\,D^4,\,\dots,\,D^r$. 

\section{Singularity classes of special multi-flags}
The Sandwich Diagram (SD in short), although completely hidden in 
the compactified definition, is a key structure underlying each special 
$m$-flag generated by, say, a rank $m+1$ distribution $D$. Much like 
in the Goursat situation, the SD allows to define an initial (rather rough) 
stratification of the germs of $D$ into so-called {\it sandwich classes}. 
(Their Goursat counterparts are {\it Kumpera-Ruiz classes} defined on 
p.\,466 in \cite{MZ1}.) Needless to say, the sandwich classes are invariant 
under the (local) automorphisms of $D$. They are further partitioned into, 
also invariant, {\it singularity classes} (defined in \cite{sin} for $m = 2$ 
and in \cite{sinbis} for arbitrary $m \ge 2$).\\

When it comes to the singularities of the special multi-flags, it is 
known for 20+ years now -- see Remark 3 in \cite{M} -- that all of them 
are realized, up to local diffeomorphisms, on the stages $P^{\,r}(\R^{m+1})$, 
$r \ge 1$, of the so-called Special Multi-Flags Monster Tower, S$m$FMT for 
short. Stages which host rank $m+1$ distributions $\D^r \subset TP^{\,r}(\R^{m+1})$ 
generating special $m$-flags of length equal to the stage number $r$.\\
So, from the singularities' point of view, without loss of generality 
one can just work with the $\D^r$'s on the stages of the S$m$FMT. 
That is to say, to classify the points of each given stage $P^{\,r}(\R^{m+1})$ 
up to the local symmetries of its proper distribution $\D^r$. (The authors 
of \cite{CH} write about `symmetries of the tower at level $k$', meaning 
diffeomorphisms of that level $k$ being at the same time symmetries of $\D^k$. 
The tower S$m$FMT goes intrinsically together with the locally universal 
distributions $\D^k$ living on its stages.)\\

\n{\it Attention.} At this point one comparison of general approaches 
to the monsters is in order. The stratifications of the multi-flags' 
monster's stages into singularity classes are -- in our opinion -- the most 
natural, canonical generalizations of the Kumpera-Ruiz stratifications in 
the Goursat Monster. Especially when one uses vector field' generators 
of the distributions $\D^r$ and visualises them at each prolongation step 
in $m+1$ affine charts of the respective projective $m$-spaces. 
They are not, however, the only possible. Alternative stratifications, 
of much different character, had been proposed in a series of works \cite{CH}, 
\cite{chs}, \cite{ccks}. The last contribution in this series somehow summarizes 
that line of constructions, by defining so-called {\it intersection loci} 
and {\it RV-strata}.\\

\n In \cite{monster} there has been undertaken an attempt at comparing the 
singularity classes and RV strata. In section 5.2 there, for $m = 2$, the union 
of classes 1.2.1.2 and 1.2.1.3 is being dissected by three RV strata. A number 
of more involved comparison examples was communicated to the authors of \cite{ccks}. 
(Already the dissection of the class 1.2.1.1.3 by RV strata is far from immediate, 
let alone that of the 1.2.1.1.1.3.) A moral to it is that singularity classes 
quickly come badly across with RV strata.\\

\n The singularity classes building up any given stage $P^{\,r}(\R^{m+1})$ 
are encoded by words \,$j_1.\,j_2\dots j_r$ over the alphabet 
$\{1,\,2,\,\dots,\,m,\,m+1\}$ such that $j_1 = 1$ and $j_{k+1} \le 1 + 
\max(j_1,\dots,\,j_k)$ for $k = 1,\dots, r-1$. (This limitation is called 
the `least upward jumps' rule', or `lujr' for short.) The approximate 
number of singularity classes in the $r$-th stage is $\frac1{(m+1)!}(m+1)^r$. 
For instance for $m = 2$ and $r = 5$ (pertinent for the present work) it is 
$41 \approx \frac1{3!}3^5$. The codimension of a class \,$j_1.\,j_2\dots j_r$ is 
$\big(\sum_{k=1}^r j_k\big) - r$. The only open dense stratum in $P^{\,r}(\R^{m+1})$ 
is 1.1$\dots$\,1 ($r$ 1's). It is homogeneous in the sense that $\D^r$ around 
every point of 1.1$\dots$\,1 is equivalent to the {\it Cartan distribution} 
(or {\it jet bundle\,} in the terminology of \cite{Y}) on the jet space $J^{\,r}(1,m)$ 
-- as observed in Theorem 1 in \cite{KR}.\\
\subsection{EKR coordinates.}
Each distribution $\D^r$ in the vicinity of a class \,$j_1.\,j_2\dots j_r$ 
is transparently visualised in its corresponding local (even: semi-global) in 
$P^{\,r}(\R^{m+1})$, polynomial coordinates, traditionally called EKR (Extended 
Kumpera-Ruiz\footnote{\,in \cite{CH} they are called Kumpera-Rubin and are $\dots$ 
not exploited to their full strength}). (Those polynomial visualisations, or 
local pseudo-normal forms, depend only on a finite number of real parameters, 
several of which can be further simplified or normalized.)
In the originating paper \cite{KR} those coordinates were introduced and used 
but not yet quite developed. Beyond the multi-flags of length 2 only so-called 
{\it truncated\,} multi-flags (related to certain systems of differential 
equations) were discussed there.\\

Those coordinates are labelled identically as the given class they visualise, 
only with {\bf bold} characters being used: $\j_1.\,\j_2\dots\j_r$. These bold 
characters stand for concrete prolongation patterns; the complete construction 
of the EKR coordinates was done in \cite{M} long before the construction of the 
singularity classes (in \cite{sin} for $m = 2$). The exact matching of the both 
concepts is far from being trivial. For $m = 2$ it is the content of Theorem 1 
in \cite{sin}; compare also Theorem 2.1 in \cite{sinbis}.\\

\n{\it Attention.} In the paper \cite{PR} certain coordinate EKR sets were constructed 
upon using intermittently only the prolongation patterns {\1} and {\2}. That is, 
only using the EKRs $\j_1.\,\j_2\dots\j_r$ with $\j_1,\,\j_2,\dots,\,\j_r \in 
\{{\1},\,{\2}\}$. Hence the main Theorem 3.2 there is incorrect (cf. in this 
respect Corollary 1 in \cite{M}).\\

The EKR coordinates, as specified for special 2-flags in section 3.3 in \cite{sin}, 
were the main technical tool in \cite{sin} and \cite{MP1}, and so they are in 
the present work, too. They are as if {\it night glasses} (Kumpera's expression) 
allowing to safely navigate in the monsters.\\

\n(Speaking uniquely about coordinate sets, coordinates equivalent to the EKRs 
were later constructed in \cite{ccks}, yet in the language of curves in the base 
and their prolongations, not of [horizontal] vector fields. So there were, from 
the one side, the EKR systems of coordinates $\j_1.\,\j_2\dots\j_r$ in $P^{\,r}(\R^{m+1})$, 
and from the other systems, denoted $C(p_1\,p_2\dots p_r)$, of coordinates in 
the $r$-th stage of the {\it Semple Tower} with $(m+1)$-dimensional base 
(= the name of S$m$FMT in algebraic geometry). The task to compare them 
stood as a genuine order of the day. A detailed discussion of these two 
equivalent ways of visualising monster' stages has been conducted 
in \cite{monster}.)\\

Reiterating already, the present text deals uniquely with special 2-flags. 
So the codes at work are the lujr words over a short alphabet $\{1,\,2,\,3\}$. 
That is, the words starting on the left with a letter 1 and such that a first 3 
(if any) goes only after a first 2 (if any).\\
The stage No 1 in S2FMT is nothing but the class 1 (cf. Theorem 1 in \cite{KR}). 
The stage No 2 is the union of the classes 1.1 (open dense) and 1.2 (of codimension 1; 
the class 1.2 is given a single local model (11) in \cite{KR}). The letter `3' shows 
up in the codes for the first time in the stage No 3, in the [name of the] class 
1.2.3 of codimension 3. 
\section{The results of the paper}
In the introductory part of \cite{ccks} the authors have written {\it `Results such 
as those in Section $5.7$ of} \,\cite{MZ2} {\it show that we can expect there to be 
infinitely many strata, that is, that there are moduli. However, we seem to be very 
far from a full understanding of where and why moduli occur.'} \,This being written 
notwithstanding an explicit example in Section 8 in \cite{MP1} of a modulus existing 
in length 7 in special $m$-flags, $m \ge 2$.\\
On the other hand, in \cite{MP2} on page 201 there was written down a dream about 
`a possible late follow-up to the work \cite{MP1}', wishfully dealing with the local 
classification problem beyond flag's length 4 completely settled in \cite{MP1}. 
In that year 2020 there were strong expectations of the finiteness of the classification 
in length 5, with possibly most difficult subtask -- handling the one-step prolongations 
of the most involved in length 4 singularity class 1.2.1.1 being itself the union 
of six orbits (cf. section 8.2 in \cite{MP2}). And there was nearly equally strong 
belief in the existence of moduli in length 6. In reality already the flags' length 5 
baffled those expectations, although the length 5 classes laden with moduli appeared 
not to be the prolongations of the class 1.2.1.1. In fact, we have recently proved 
\begin{thm}\label{main}
The local classification of special \,$2$-flags in length \,$5$ is not finite. 
In fact, in each of the singularity classes $1.2.1.2.1$, $1.2.2.1.2$ and $1.2.3.1.2$ 
there resides a single numerical modulus of the local classification. The moduli 
in the classes $1.2.1.2.1$ and $1.2.2.1.2$ live in codimension \,$3$, and in 
the class $1.2.3.1.2$ -- in codimension \,$4$.\\
The local classification within the remaining \,$38$ singularity classes 
existing in length \,$5$ is finite.
\end{thm}
The proof of this theorem is spread over the three consecutive sections here 
below, each section devoted to one of the classes put in relief in the theorem's 
wording. (In Appendix we just comment on theorem's last statement.)\\

\n The formal skeletons of proofs in the three sections are similar. 
This particularly applies to sections \ref{12312} and \ref{12212} because 
the codes of the respective singularity classes differ by just one letter. 
Symbols used in the proofs in these three sections do overlap. 
The same letters stand in them for the EKR coordinates actually in use, but remember 
that those coordinate sets critically depend -- according to Theorem 1 in \cite{sin} 
-- on the relevant singularity classes. Also the letters for the components of 
conjugating diffeomorphisms are the same in the three sections (and always 
closely follow {\it graphically} the respective variables).\\

\n Notwithstanding all this, proofs in these sections differ substantially in 
their essences and deserve, in our opinion, to be included all three in this paper. 
Of particular interest may be perceived the modulus found in the singularity class 
1.2.1.2.1. This finding is to be compared with the first ever finding of a modulus 
in special 2-flags, still in \cite{MP1}, in a distant grandson 1.2.1.2.1.2.1 of 
our actual 1.2.1.2.1.\\

As for the last statement in Theorem \ref{main}, its justification follows from 
the many partial outcomes of computations --- always made with the use of \cite{W} 
and with the initial classification-related information in length 4 taken from 
\cite{MP1} --- of the infinitesimal symmetries of the germs of $\D^5$ in the 
remaining singularity classes in length 5 not listed in the theorem.\\
To illustrate that kind of argument, in Appendix we deal with one chosen class 
1.2.1.3.1 and give the arguments that it consists of a finite number of orbits 
of the local classification by using both a)\,the general form of an infinitesimal 
symmetry of [a part of] 1.2.1.3.1 evaluated at $0 \in \R^{13}$, and b)\,explicit 
rescalings of EKR variables normalizing pairs of new incoming parameters.\\

Prior to all concrete arguments there follow short historical and introductory 
remarks.\\

An important step towards clarifying the (local) status of 2-flags' lengths 5 and 6 
were the recursive formulae for the {\it infinitesimal symmetries} of the special 
2-flags put forward in \cite{MP2}. Yet a real breakthrough occurred only during 
the winter 2023/24 when A. Weber effectively implemented in the {\it Wolfram 
Mathematica}, \cite{W}, the formulae from \cite{MP2}.\\
When it came to the application of that implementation to the 2-flags of length 5, 
the exact local models in length 4, known since long, played the role of the initial 
conditions. Then the components of the i.s.'s at level 5 were quickly given explicit 
algebraic expressions. (Compare central in this respect Lemmas 2, 3 and 4 in \cite{MP2}. 
Reiterating, a closer glimpse of the technique based on the computation of all 
the i.s.'s of a given distribution, here $\D^5$, is given in the paper's 
end in Appendix.)\\

The output of Weber's computations was the list of {\it three} singularity 
classes, each {\it likely} hosting a single modulus of the local classification. 
Meaning that the newly incoming parameters in the respective EKR's could not be 
moved, keeping all the previous parameter values frozen, by means of diffeomorphisms 
(automorphisms of $\D^5$) embeddable in the flows of autonomous vector fields. 
Speaking not, however, about possible automorphisms of $\D^5$ {\it not\,} 
embeddable into such flows that also ought to be considered. The classes 
identified by Weber were 
\begin{itemize}
\item 1.2.1.2.1,\\
\item 1.2.2.1.2,\\
\item 1.2.3.1.2
\end{itemize}
(reiterating, the total of singularity classes of special 2-flags 
in length 5 is 41, according to Proposition 3 in \cite{sin}). 
\subsection{A modulus in the class 1.2.3.1.2\,.}\label{12312}
On the above short list, the shortest line of arguments happens to concern 
the singularity class 1.2.3.1.2. In this case the `father' class of length 4 
is 1.2.3.1 which is known (cf. pages 24--25 in \cite{MP1}) to be built up 
by four orbits. Moreover, Weber's computations \cite{W} prompted even more 
precisely that a plausible modulus in 1.2.3.1.2 could only reside in the 
one-step prolongation of the thickest first orbit inside 1.2.3.1, called 
in \cite{MP1} `transverse generic' (the lion's part of the entire father 
class 1.2.3.1).\\

Before all technicalities, agree, for the purpose of the proof in this 
and the forthcoming sections, not to watch (in this section in the glasses 
{\1}.{\2}.{\3}.{\1}.{\2}) the universal distribution $\D^5$ but instead 
just to consider two distribution germs $\D_c^5$ and $\D_{\t c}^5$, 
where $\D_a^5$ ($a \in \R$) is the distribution 
\begin{align*}\label{312}
dx0 - x1dt        &= 0   &  dy0 - y1dt       &= 0\\
dt - x2\,dx1      &= 0   &  dy1 - y2\,dx1    &= 0\\
dx1 -  x3\,dy2    &= 0   &  dx2 - y3dy2      &= 0\\
dx3 - (1+x4)dy2   &= 0   &  dy3 - (1+y4)dy2  &= 0\\
dy2 - x5\,dx4     &= 0   &  dy4 - (a+y5)dx4  &= 0\,. 
\end{align*}
viewed as the germ at $0 \in \R^{13}\big(t,\,x0,\,y0,\dots,\,x5,\,y5\big)$. 
(This is its description in the form of a set of 10 Pfaffian equations, in 
the style of the pioneering work \cite{KR} and of \cite{MP1}.) In explanation, 
we have taken here a duly normalized EKR {\1}.{\2}.{\3}.{\1} representative 
of the mentioned transverse generic orbit and extended it by a single pair 
of new Pfaffian equations, remembering that the last letter in the code 
of the `son' class 1.2.3.1.2 is `2' (cf. in \cite{sin} Theorem 1 coupled 
with section 3.3 there).\\
However, the vector field generators' description of each distribution 
$\D_a^5$ is handier. Skipping here and in the sequel the symbol `span': 
\be\label{zzz}
\D_a^5 = \big(Z_{1,a}\,,\,Z_2,\,Z_3\big), 
\ee
where (we systematically use a shorthand notation of the type $\p_t = \p/\p t$ etc.)
\begin{align}
Z_{1,a} = &\ x5(x3\big(x2\Big(\p_t + x1\,\p_{x0} + y1\,\p_{y0}\Big) + \p_{x1} + y2\,\p_{y1}\big)\nonumber\\
&+ y3\,\p_{x2} + \p_{y3} + (1 + x4)\p_{x3} + (1 + y4)\p_{y3}) + \p_{x4} + (a + y5)\p_{y4}\,,\nonumber\\
&\qquad\qquad Z_2 = \p_{x5}\,,\qquad\qquad Z_3 = \p_{y5}\,.
\end{align}
The aim is to show that for every two {\bf different} values $a = c$ and $a = \t c$ 
the distribution germs \eqref{zzz} are {\bf not equivalent}. To this end we work with 
an arbitrary local diffeomorphism in $P^5(\R^3)$ written in the EKR coordinates 
\[
\varPhi \;= \;\big(T,\,X0,\,Y0,\,X1,\,Y1,\dots,\,X5,\,Y5\big):
\;\big(\R^{13},\,0\big) \hookleftarrow
\]
and supposed to conjugate these two local distributions: 
\be\label{key}
d\,\varPhi(p)\D_c^5(p) = \D_{\t c}^5\big(\varPhi(p)\big)
\ee
for all points $p$ close to $0 \in \R^{13}$. 
Our goal is to arrive at the conclusion $c = \t c$.\\

\n From the general considerations in \cite{sin} it is well known that 
\begin{itemize}
\item the functions $T,\,X0,\,Y0$ depend only on $t,\,x0,\,y0$\,,
\item for $1 \le j \le 5$, functions $Xj,\,Y\!j$ depend 
only on $t,\,x0,\,y0,\,x1,\,y1,\dots,\,xj,\,yj$.
\end{itemize}
Moreover, in the discussed situation the 2nd, 3rd and 5th entries 
in the relevant EKR code are not {\1} (the class 1.2.3.1.2 sits 
in the bigger sandwich class 1.\u{2}.\u{2}.1.\u{2} as defined in 
section 3.2 in \cite{sin}). Hence much more is known -- see Appendix A 
in \cite{sin} -- about the respective components $X2,\,\,X3$, and $X5$ 
of the diffeomorphism $\varPhi$. Namely, $\varPhi$ has to preserve the loci 
$\{x2 = 0\}$, $\{x3 = 0\}$ and $\{x5 = 0\}$ of the pointwise inclusions 
in the respective sandwiches No 2, 3 and 5. That is to say 
\begin{itemize}
\item $X2(t,\,x0,\,y0,\,x1,\,y1,\,x2,\,y2) = 
x2\,K(t,\,x0,\,y0,\,x1,\,y1,\,x2,\,y2)$\,,
\item $X3(t,\,x0,\dots,\,y3) = x3\,H(t,\,x0,\dots,\,x3,\,y3)$\,,
\item $X5(t,\,x0,\dots,\,y5) = x5\,G(t,\,x0,\dots,\,x5,\,y5)$
\end{itemize}
for certain functions $G,\,H,\,K$. Observe -- this is important -- that 
all these three functions $G,\,H,\,K$ are invertible at 0. Moreover still, 
the preservation of $\D^5$ by $\varPhi$ clearly implies the preservation 
by $\varPhi$ of the Cauchy characteristic distribution 
\be\label{Cauchy}
L\Big(\big[\D^5,\,\D^5\big]\Big) = \big(Z_2,\,Z_3\big) = \big(\p_{x5},\,\p_{y5}\big)
\ee
which is common for both local structures $\D_c^5$ and $\D_{\t c}^5$. 
Therefore, along with the conjugating diffeomorphism $\varPhi$ there must exist 
an invertible at 0 function germ $f$, $f\0 \;\ne 0$, that, altogether, 
\[
d\,\varPhi(p)Z_{1,c}(p) = f(p)Z_{1,\t c}\big(\varPhi(p)\big) + 
(\ast)Z_2\big(\varPhi(p)\big) + (\ast\ast)Z_3\big(\varPhi(p)\big)
\] 
for all $p$'s close to $0 \in \R^{13}$. The same equation in an expanded form reads
\be\label{mo1}
d\,\varPhi(p)
\ba{r}x5\!\left(\!\!\!\!
       \ba{r}
       x3\!\left(\!\!\!\ba{r}
       \left.x2\!\left(\!\!\ba{c}1\\x1\\y1\ea\right.\right]\mbox{\hskip.0005mm}\\
       \left.\ba{c}1\\y2\ea\mbox{\hskip1mm}\right]\mbox{\hskip.001mm}
       \mbox{\hskip.001mm}\ea\right.\\
\left.\ba{c}y3\mbox{\hskip2mm}\\1\mbox{\hskip2mm}\ea\right]\mbox{\hskip2mm}\\
\left.\ba{r}1 + x4\mbox{\hskip2mm}\\1 + y4\mbox{\hskip2mm}\ea\right]
\mbox{\hskip2mm}\ea\right.\\
\left.\ba{c}1\mbox{\hskip2.6mm}\\c + y5\mbox{\hskip3mm}\ea\right]
\mbox{\hskip4.2mm}\\
\left.\ba{c}0\mbox{\hskip4.1mm}\\0\mbox{\hskip4.1mm}\ea\right]
\mbox{\hskip4.2mm}\ea
=\quad 
\ba{r}f(p)\!\left(\!\!\!\!\ba{r}
      x5\,G\!\left(\!\!\!\!
       \ba{r}
       x3\,H\!\left(\!\!\!\!\ba{r}
             \left.x2\,K\!\left(\!\!\ba{c}1\\X1\\Y1\ea\right.\right]\\
           \left.\ba{c}1\\Y2\ea\mbox{\hskip1mm}\right]\mbox{\hskip.005mm}\\
                           \ea\right.\\
\left.\ba{c}Y3\mbox{\hskip2mm}\\1\mbox{\hskip2mm}\ea\right]\mbox{\hskip2mm}\\
\left.\ba{r}1 + X4\mbox{\hskip2mm}\\1 + Y4\mbox{\hskip2mm}\ea\right]
\mbox{\hskip2mm}\ea\right.\\  
\left.\ba{c}1\mbox{\hskip2mm}\\\t{c} + Y5\mbox{\hskip2.4mm}\ea\right]
\mbox{\hskip4.2mm}\\
\ea\right.\\
\left.\ba{c}(\ast)\mbox{\hskip4.2mm}\\ (\ast\ast)\mbox{\hskip5mm}\ea\right]
\mbox{\hskip6.5mm}
\ea
\ee
where $p \,= (t,\,x0,\,y0,\dots,\,x5,\,y5)$ and, for bigger transparence, the 
arguments in the functions $G,\,H,\,K,\,X1,\dots,\,Y5$ on the RHS of \eqref{mo1} 
are not written. This huge vector relation entails the set of 11 scalar equations 
equalling the coefficients, on the both sides of \eqref{mo1}, at the components 
$\p_t$, $\p_{x0},\,\dots$, $\p_{x4}$, $\p_{y4}$. 
(We disregard the last two components in \eqref{mo1} -- the components 
in the directions spanning the Cauchy characteristics' subdistribution 
$L\Big(\big[\D^5,\,\D^5\big]\Big)$, cf. \eqref{Cauchy}.)\\

\n In view of the first NINE components of the diffeomorphism $\varPhi$ depending 
only, altogether, on $t,\,x0,\dots,\,x3,\,y3$, the upper NINE among these scalar 
equations can be divided sidewise by $x5$. Likewise and additionally, the upper 
FIVE among them can be divided by $x3$, and the first THREE equations -- 
additionally divided by $x2$.\\
Agree to call thus simplified equations `level $T$', `level $X0$', `level $X4$', 
etc, in function of the row of the matrix $d\,\varPhi(p)$ being involved. 
For instance, the level $T$ equation is the equalling of the coefficients 
on the both sides of \eqref{mo1} at the $\p_t$--component, divided sidewise 
by the product of variables $x2\,x3\,x5$. 
\vskip2.5mm
\n Ad rem now, from level $X4$ we infer that 
\be\label{13}
X4_{\,x4} + c\,X4_{\,y4}\0 \,\,= f\0\,, 
\ee
and from level $Y4$ we know 
\be\label{14}
Y4_{\,x4} + c\,Y4_{\,y4}\0 \,\,= \,{\t c}\,f\0\,.
\ee
Level $X3$ yields 
\be\label{15}
(\ast)x3 + (H + x3\,H_{\,x3})(1 + x4) + x3\,H_{\,y3}(1 + y4) = fG(1 + X4), 
\ee
hence 
\be\label{16}
H\0 \,\,= fG\0\,.
\ee
The function $fG$, got from level $Y2$, depends only on $t,\,x0,\,y0,\dots,\,x3,\,y3$. 
Hence, differentiating \eqref{15} sidewise with respect to $x4$ at 0 and then using 
\eqref{16}, one gets 
\be\label{17}
X4_{\,x4}\0 \,\,= 1. 
\ee
While differentiating \eqref{15} sidewise with respect to $y4$ at 0, 
\be\label{18}
X4_{\,y4}\0 \,\,= 0\,.
\ee
Joining the pieces \eqref{13}, \eqref{17} and \eqref{18}, 
\be\label{19}
f\0 \,\,= 1\,.
\ee
Now it is the identity at level $Y3$ that comes in handy: 
\be\label{20}
(\ast)\,x3 + (\ast)\,y3 +  Y3_{\,y2} + Y3_{\,x3}(1 + x4) + Y3_{\,y3}(1 + y4) = fG(1 + Y\!4)\,.
\ee
Remembering that $fG$ does not depend on $x4$ and differentiating \eqref{20} sidewise 
with respect to $x4$ at 0, 
\[
fG\,\,Y\!4_{\,x4}\0 \,\,= \,\big(fG(1 + Y\!4)\big)_{x4}\0 \,\,= Y3_{\,x3}\0 \,\,= \,0
\]
because level $X2$ says that 
\be\label{20bis}
{\rm the}\ \ {\rm function}\ \ fG\,\,Y3\ \ {\rm is}\ \ {\rm a}\ \ {\rm functional}
\ \ {\rm combination}\ \ {\rm of}\ \ x2\ \ {\rm and}\ \ y3\,,
\ee
whence $fG\,\,Y3_{\,x3}\0 \,\,= \,(fG\,\,Y3)_{x3}\0 \,\,= \,0$. 
Therefore -- see also \eqref{19} -- the key relation \eqref{14} reduces to 
\be\label{21}
c\,\,Y\!4_{\,y4}\0 \,\,= \,{\t c}.
\ee
But \eqref{20} can also be differentiated sidewise with respect to $y4$ at 0: 
\be\label{22}
Y3_{\,y3}\0 \,\,= \,fG\,\,Y\!4_{\,y4}\0. 
\ee
Now the end of proof is virtually around the corner, because \eqref{20} means at 0 
\[
Y3_{\,y2} + Y3_{\,x3} + Y3_{\,y3}\0 \,\,= \,fG\0
\]
and the first two summands on the left hand side vanish because of \eqref{20bis}. 
That is, 
\be\label{23}
Y3_{\,y3}\0 \,\,= fG\0 \,\,\ne \,0\,.
\ee
Now we are done, because the relations \eqref{22} and \eqref{23} 
taken together imply $Y\!4_{\,y4}\0 \,\,= \,1$, or, 
via \eqref{21}, $c = \t c$.\\

\n A modulus in the singularity class 1.2.3.1.2 has now been justified.\\

\n{\bf Remark.} A slightly shorter proof would overtly use the fact that 
the third letter in the code 1.2.3.1.2 is `3'. It implies that the component 
function $Y3$ of the conjugating diffeomorphism $\varPhi$ is a functional combination 
of $x3$ and $y3$: $Y3 = x3\,M + y3\,N$, $N(0) \ne 0$, cf. Appendix A in \cite{sin}. 
However, the fact $j_3 = 3$ {\it has} been used in the above proof: the precious 
information \eqref{20bis} is due to it. In the proof for the class 1.2.2.1.2 
in the section that follow the knowledge about the respective function $Y3$ 
will be incomparably scarcer.
\subsection{A modulus in the class 1.2.2.1.2\,.}\label{12212}
The first point in this part of the proof is, not surprisingly, an inspection 
of the `father' class 1.2.2.1. It also consists of four orbits (cf. pages 22--24 
in \cite{MP1}). Much like in the previous section, also here Weber's computations 
\cite{W} prompted that a plausible modulus in 1.2.2.1.2 could only reside in the 
prolongation of the thickest first orbit inside 1.2.2.1 -- the lion's part of the 
entire father class 1.2.2.1. So we take the normalized EKR normal form of that 
orbit and extend it by a pair of Pfaffian equations written in concordance 
with the fifth letter `2' in the code 1.2.2.1.2: 
\begin{align}\label{212}
dx0 - x1dt         &= 0   &  dy0 - y1dt       &= 0\nonumber\\
dt   - x2\,dx1     &= 0   &  dy1 - y2\,dx1    &= 0\nonumber\\
dx1 - x3\,dx2      &= 0   &  dy2 - y3\,dx2    &= 0\nonumber\\
dx3 - (1+x4)\,dx2  &= 0   &  dy3 - (1+y4)\,dx2&= 0\\
dx2 - x5\,dx4      &= 0   &  dy4 - (c+y5)\,dx4&= 0\,.\nonumber
\end{align}
Here $c \in \R$ is a real parameter and these objects are understood 
as germs at $0 \in \R^{13}(t,\,x0,\,y0,\,x1,\,y1,\dots,\,x5,\,y5)$. 
We want to show that for every two different values $c$ and $\t c$ 
{\bf the distribution germs \eqref{212} are non-equivalent}. 
To this end we analyze an arbitrary local diffeomorphism 
\[
\varPhi \;= \;\big(T,\,X0,\,Y0,\,X1,\,Y1,\dots,\,X5,\,Y5\big):\;\big(\R^{13},\,0\big) 
\hookleftarrow
\]
sending the distribution featuring a constant $c$ 
to the one featuring a constant $\t c$.\\

The general type information about such a diffeo $\varPhi$ is the same as in section 
\ref{12312} because both singularity classes 1.2.2.1.2 and 1.2.3.1.2 are included in 
the same sandwich class 1.\u{2}.\u{2}.1.\u{2} which is solely responsible for the particular 
form of the components $X2,\,X3,\,X5$ of $\varPhi$ put forward in section \ref{12312}. 
However, the basic vector equation \eqref{mo2}implied by the conjugation effectued 
by the diffeomorphism $\varPhi$ is different from \eqref{mo1}, simply because the EKR 
pseudo-normal forms are now different -- and they faithfully shadow the respective 
singularity classes (Theorem 1 in \cite{sin}). Here is that actual vector equation 
laden with information sufficient to conclude the proof of Theorem \ref{main} for 
the class 1.2.2.1.2\,. 
\be\label{mo2}
d\,\varPhi(p)
\ba{r}x5\!\left(\!\!\!\!
       \ba{r}
       x3\!\left(\!\!\!\ba{r}
       \left.x2\!\left(\!\!\ba{c}1\\x1\\y1\ea\right.\right]\mbox{\hskip.0005mm}\\
       \left.\ba{c}1\\y2\ea\mbox{\hskip1mm}\right]\mbox{\hskip.001mm}
       \mbox{\hskip.001mm}\ea\right.\\
\left.\ba{c}1\mbox{\hskip2mm}\\y3\mbox{\hskip2mm}\ea\right]\mbox{\hskip2mm}\\
\left.\ba{r}1 + x4\mbox{\hskip2mm}\\1 + y4\mbox{\hskip2mm}\ea\right]
\mbox{\hskip2mm}\ea\right.\\
\left.\ba{c}1\mbox{\hskip2.6mm}\\c + y5\mbox{\hskip3mm}\ea\right]
\mbox{\hskip4.2mm}\\
\left.\ba{c}0\mbox{\hskip4.1mm}\\0\mbox{\hskip4.1mm}\ea\right]
\mbox{\hskip4.2mm}\ea
=\quad 
\ba{r}f(p)\!\left(\!\!\!\!\ba{r}
      x5\,G\!\left(\!\!\!\!
       \ba{r}
       x3\,H\!\left(\!\!\!\!\ba{r}
             \left.x2\,K\!\left(\!\!\ba{c}1\\X1\\Y1\ea\right.\right]\\
           \left.\ba{c}1\\Y2\ea\mbox{\hskip1mm}\right]\mbox{\hskip.005mm}\\
                           \ea\right.\\
\left.\ba{c}1\mbox{\hskip2mm}\\Y3\mbox{\hskip2mm}\ea\right]\mbox{\hskip2mm}\\
\left.\ba{r}1 + X4\mbox{\hskip2mm}\\1 + Y4\mbox{\hskip2mm}\ea\right]
\mbox{\hskip2mm}\ea\right.\\  
\left.\ba{c}1\mbox{\hskip2mm}\\\t{c} + Y5\mbox{\hskip2.4mm}\ea\right]
\mbox{\hskip4.2mm}\\
\ea\right.\\
\left.\ba{c}*\mbox{\hskip4.2mm}\\ \ast\mbox{\hskip5mm}\ea\right]
\mbox{\hskip6.5mm}
\ea
\ee
Also the distinguished subsets of scalar equations implied by \eqref{mo2} 
and the names attached to particular equations are formally identical to 
those in the previous section, the same underlying sandwich class 
1.\u{2}.\u{2}.1.\u{2} being responsible for it.\\
Moreover, the local conclusions \eqref{13} through \eqref{19} in the previous 
section carry over verbatim for the components of the actual diffeomorphism 
$\varPhi$ because the last two letters in the classes 1.2.3.1.2 and 1.2.2.1.2 
are identical.\\

\n Now it is the identity at level $Y3$ that comes in handy: 
\be\label{31}
(\ast)x3 + Y3_{\,x2} + (\ast)y3 + Y3_{\,x3}(1 + x4) + Y3_{\,y3}(1 + y4) = fG(1 + Y4)\,.
\ee
In particular 
\[
Y3_{\,x2} + Y3_{\,x3} + Y3_{\,y3}\0 \,\,= fG\0\,,
\]
or else 
\be\label{32}
\frac1{fG}\Big(Y3_{\,x2} + Y3_{\,x3} + Y3_{\,y3}\Big)\0\,\,= 1\,.
\ee
Upon differentiating \eqref{31} sidewise at 0: firstly with respect to $x4$, 
then with respect to $y4$, one gets 
\[
Y4_{\,x4}\0 \,\,= (fG)^{-1}Y3_{\,x3}\0 \qquad\quad {\rm and}\qquad\quad 
Y4_{\,y4}\0 \,\,= (fG)^{-1}Y3_{\,y3}\0\,.
\]
These identities together with \eqref{19} transform relation \eqref{14} into 
\be\label{33}
\frac1{fG}\Big(Y3_{\,x3} + c\,Y3_{\,y3}\Big)\0\,\,= \,\t c\,. 
\ee
{\bf Remark.} The couple of relations \eqref{32} -- \eqref{33} is as if 
a `chuck' that will inevitably lead to the conclusion $c = \t c$.\\

\n It prompts by itself to replace $Y3$ in \eqref{32} -- \eqref{33} by $Y2$, 
and level $Y2$ reads, handy enough, 
\be\label{34}
x3\big(x2(\cdots) + Y2_{\,x1} + y2(\cdots)\big) + Y2_{\,x2} + y3\,Y2_{\,y2} = fG\,Y3\,.
\ee
It remains to differentiate \eqref{34} sidewise, always at 0, with respect to $x2$, 
then to $x3$, and then to $y3$, obtaining respectively: 
\[
Y2_{\,x2x2}\0 \,\,= \,fG\,\,Y3_{\,x2}\0\,,
\]
\[
Y2_{\,x1}\0 \,\,= \,fG\,\,Y3_{\,x3}\0\,,
\]
\[
Y2_{\,y2}\0 \,\,= \,fG\,\,Y3_{\,y3}\0\,.
\]
With these at hand, \eqref{32} becomes 
\be\label{35}
\Big(\frac1{fG}\Big)^2\Big(Y2_{\,x2x2} + Y2_{\,x1} + Y2_{\,y2}\Big)\0\,\,\,= \,1\,, 
\ee
and \eqref{33} becomes 
\be\label{36}
\Big(\frac1{fG}\Big)^2\Big(Y2_{\,x1} + c\,Y2_{\,y2}\Big)\0\,\,\,= \,\t c\,.
\ee
\vskip8mm
\n{\bf Observation.} (a) \,$Y2_{\,x2x2}\0 \,\,= \,0$\,,\\
(b) \,\,$Y2_{\,x1}\0 \,\,= \,0$\,.\\

\n Proof of (a).\\
Because $Y2 = \frac{fGH\,Y2}{fGH}$, one uses: level $Y1$ in \eqref{mo2} for the numerator 
and level $X1$ for the denominator. In fact, for computing the second derivative $Y2_{\,x2x2}\0$\,, 
\[
Y2 \equiv \frac{x2\,Y1_t + Y1_{\,x1}}{x2\,X1_t + X1_{\,x1}}\qquad {\rm and}\qquad 
Y2_{\,x2} \equiv \frac{Y1_t\,X1_{\,x1} - Y1_{\,x1}\,X1_t}{\big(x2\,X1_t + X1_{\,x1}\big)^2}\,.
\]
The numerator in the second equivalence above vanishes at 0 since so does the entire 
fraction evaluated at 0, $Y2_{\,x2}\0 \,\,= \, fG\,Y3\0\,\,= \,0$. And that numerator 
does not depend whatsoever on $x2$. This proves already (a).\\

\n Proof of (b). 
\[
Y2_{\,x1}\0 \,\,= \,\left(\frac{Y1_{\,x1}}{X1_{\,x1}}\right)_{x1}\0\,\,= 
\,\frac{Y1_{\,x1x1}X1_{\,x1} - Y1_{\,x1}X1_{\,x1x1}}{\big(X1_{\,x1}\big)^2}\0\,.
\]
Both summands in the numerator on the right-most side vanish. The subtrahend, because 
\be\label{37}
Y1_{\,x1}\0\,\,= \,fGH\,Y2\0\,\,= \,0\,. 
\ee
The minuend, because 
\be\label{38}
Y1_{\,x1x1}\0 \,\,= \,\left(\frac{Y0_t + x1Y0_{\,x0}}{T_t + x1T_{x0}}\right)_{x1x1}\0 \,\,= 
\,\left(\frac{T_tY0_{\,x0} - T_{x0}Y0_{\,t}}{\big(T_t + x1T_{x0}\big)^2}\right)_{x1}\0 \,\,= \,0\,.
\ee
This last equality necessitates an explanation. Namely, $T_tY0_{\,x0} - T_{x0}Y0_{\,t}$ 
does not depend on $x1$ and vanishes at 0. Indeed, its independence of $x1$ is visible and 
$Y0_{\,t}\0 \,\,= \,0$ because level $Y0$ in \eqref{mo2} shows that it is $fGHK\,Y1\0$. 
Endly, $Y0_{\,x0}\0 \,\,= \,0$ because, after \eqref{37}, 
\[
0 = \,Y1_{\,x1}\0 \,\,= \,\frac{T_tY0_{\,x0} - T_{x0}Y0_{\,t}}{\big(T_t + x1T_{x0}\big)^2}\0
\]
from the course of computation in the line \eqref{38}, and the equality $Y0_{\,t}\0 \,\,= \,0$ 
has been ascertained a moment ago. Hence 
\[
0 = \frac{Y0_{\,x0}}{T_t}\0\,.
\]
This ends the proof of (b) and of the entire Observation.\\

\n In view of Observation, the equality \eqref{35} boils down to 
\be\label{39}
\Big(\frac1{fG}\Big)^2\,Y2_{\,y2}\0\,\,= 1\,, 
\ee
while \eqref{36} boils down to 
\be\label{40}
\Big(\frac1{fG}\Big)^2\,c\,\,Y2_{\,y2}\0\,\,\,= \,\t c\,.
\ee
Now \eqref{39} and \eqref{40} yield $c = \t c$.\\
A modulus in the singularity class 1.2.2.1.2 has been justified. 
\subsection{A modulus in the class 1.2.1.2.1\,.}\label{12121}
As in the previous parts of the proof of Theorem \ref{main}, an inspection 
of the father class 1.2.1.2 is obligatory. This time Weber's hint, \cite{W}, 
has been to focus (one more time!) on the lion's part orbit in 1.2.1.2 having 
the normal form description (37) on page 28 in \cite{MP1}. 
However, this has not yet fully determined the needed pseudo-normal forms 
to work with inside 1.2.1.2.1, because now the last code's letter is a general 
position' letter `1'. A super hint from Weber's infinitesimal symmetries' computations 
was to not concentrate now on the general points prolonging the orbit prompted at 
level 4,\footnote{\,themselves also building up just one orbit inside 1.2.1.2.1} 
but instead to deal with points where $\D^5$ is {\it tangent\,} to the locus 
of the father 1.2.1.2 singularity. This means, in the relevant EKR glasses, 
the absence of an additive constant standing next to the variable $x5$.
In the down-to-earth practice it means that we are to work with the 
following 1-parameter subfamily of EKR's {\1}.{\2}.{\1}.{\2}.{\1}\,:
\begin{align}\label{121}
dx0 - x1dt         &= 0   &  dy0 - y1dt       &= 0\nonumber\\
dt   - x2\,dx1     &= 0   &  dy1 - y2\,dx1    &= 0\nonumber\\
dx2 - (1+x3)\,dx1  &= 0   &  dy2 - y3\,dx1    &= 0\nonumber\\
dx1 - x4\,dx3      &= 0   &  dy3 - (1+y4)\,dx3&= 0\\
dx4 - x5\,dx3      &= 0   &  dy4 - (c+y5)\,dx3&= 0\,.\nonumber
\end{align}
Geometrically this is a second order singularity: the mentioned 
tangency condition adds one codimension to the codimension two 
of the class 1.2.1.2.1: $2 + 1 = \,3$.\\
Here $c \in \R$ is an arbitrary real parameter and these objects are 
understood as germs at $0 \in \R^{13}(t,\,x0,\,y0,\,x1,\,y1,\dots,\,x5,\,y5)$. 
We want to show that for every two different values $c$ and $\t c$ 
the distribution germs \eqref{121} are non-equivalent. To show this 
we start to analyze an arbitrary local diffeomorphism 
\[
\varPhi \;= \;\big(T,\,X0,\,Y0,\,X1,\,Y1,\dots,\,X5,\,Y5\big):\;\big(\R^{13},\,0\big) 
\hookleftarrow
\]
assumed to conjugate these two distributions. The general limitations on the function 
components of such a $\Phi$ are exactly as in the two previous sections. As to the 
additional limitations, in the discussed situation only the 2nd and 4th letters 
in the code 1.2.1.2.1 are not `1' -- the class sits inside the sandwich class 
1.\u{2}.1.\u{2}.1\,. Hence it is known from the beginning about the components 
$X2$ and $X4$ that 
\begin{itemize}
\item $X2(t,\,x0,\,y0,\,x1,\,y1,\,x2,\,y2) = 
x2\,H(t,\,x0,\,y0,\,x1,\,y1,\,x2,\,y2)$\,,
\item $X4(t,\,x0,\dots,\,y4) = x4\,G(t,\,x0,\dots,\,x4,\,y4)$
\end{itemize}
for certain invertible at 0 functions $G,\,H$. Moreover, identically as in 
the two previous sections, there must exist an invertible at 0 function $f$, 
$f\0 \;\ne 0$, such that, altogether, 
\be\label{mo3}
d\,\varPhi(p)
\ba{r}\!\!\!\!
       \ba{r}
       x4\!\left(\!\!\!\ba{r}
       \left.x2\!\left(\!\!\ba{c}1\\x1\\y1\ea\right.\right]\mbox{\hskip.0005mm}\\
       \left.\ba{c}1\\y2\ea\mbox{\hskip1mm}\right]\mbox{\hskip.001mm}
       \mbox{\hskip.001mm}\\
\left.\ba{c}1 + x3\mbox{\hskip2mm}\\y3\mbox{\hskip2mm}\ea\right]\ea\right.\\
\left.\ba{r}1\mbox{\hskip2mm}\\1 + y4\mbox{\hskip2mm}\ea\right]\mbox{\hskip2mm}\ea\\
\left.\ba{c}x5\mbox{\hskip2.6mm}\\c + y5\mbox{\hskip3mm}\ea\right]
\mbox{\hskip4.2mm}\\
\left.\ba{c}0\mbox{\hskip4.1mm}\\0\mbox{\hskip4.1mm}\ea\right]
\mbox{\hskip4.2mm}\ea
=\quad 
\ba{r}f(p)\!\left(\!\!\!\!\ba{r}
      \!\!\!\!
       \ba{r}
       \,\,x4\,G\!\left(\!\!\!\!\ba{r}         
			 \left.x2\,H\!\left(\!\!\ba{c}1\\X1\\Y1\ea\right.\right]\\
			 \left.\ba{c}1\\Y2\ea\mbox{\hskip1mm}\right]\mbox{\hskip.005mm}\\
       \left.\ba{c}1 + X3\mbox{\hskip2mm}\\Y3\ea\right]
			 \ea\right.\\
\left.\ba{r}1\mbox{\hskip2mm}\\1 + Y4\mbox{\hskip2mm}\ea\right]
\mbox{\hskip2mm}\ea\\  
\left.\ba{c}X5\mbox{\hskip2mm}\\\t{c} + Y5\mbox{\hskip2.4mm}\ea\right]
\mbox{\hskip4.2mm}\\
\ea\right.\\
\left.\ba{c}*\mbox{\hskip4.2mm}\\ \ast\mbox{\hskip5mm}\ea\right]
\mbox{\hskip6.5mm}
\ea
\ee
where $p \,= (t,\,x0,\,y0,\dots,\,x5,\,y5)$. In view of the first SEVEN components 
of the diffeomorphism $\varPhi$ depending only on $t,\,x0,\dots,\,x2,\,y2$, the upper 
SEVEN among the scalar equations contained (or: hidden) in \eqref{mo3} can be divided 
sidewise by $x4$. While the upper THREE among them can be divided sidewise by 
the product of variables $x2\,x4$.\\
Agree, as in the preceding sections, to call thus simplified scalar equations 
`level $T$', `level $X0$', `level $X4$', etc, in function of the row of $d\,\varPhi(p)$ 
being involved. For instance, the level $Y0$ equation is the $\p_{y0}$--component 
scalar equation in \eqref{mo3} divided sidewise by $x2\,x4$.\\

\n The relation binding the constants $c$ and $\t c$ is encoded in level $Y4$: 
\be\label{41}
Y4_{\,x3} + Y4_{\,y3} + c\,Y4_{\,y4}\0 \,\,= \,{\t c}\,f\0. 
\ee
{\bf Observation 1.} 
\[
Y4_{\,x3}\0 \,\,\,= \,Y4_{\,y3}\0 \,\,\,= \,0\,.
\]
Proof. The component function $Y4$ is expressed in terms of the function $Y3$ 
in level $Y3$: 
\be\label{42}
x4(\ast) + Y3_{\,x3} + Y3_{\,y3}(1 + y4) = f(1 + Y4).
\ee
In \eqref{42} there also show up the functional coefficient $f$. 
It is explicitly got in level $X3$ in \eqref{mo3}: 
\be\label{43}
f \,= \,x4(\ast) + X3_{\,x3} + X3_{\,y3}(1 + y4).
\ee
Level $X1$ in \eqref{mo3} tells us that the function coefficient 
\be\label{28}
fG\ {\rm depends}\ {\rm only}\ {\rm on}\ t,\,x0,\,y0,\dots,\,x2,\,y2.
\ee
In consequence, both the component functions $X3$ and $Y3$ showing up in levels 
$X2$ and $Y2$, respectively, are {\it affine} in the variable $x3$. Knowing this, 
we differentiate sidewise with respect to $x3$ at 0 the relations \eqref{42} 
and \eqref{43}: 
\[
f\,Y4_{\,x3}\0 \,\,= \,Y3_{\,x3x3} + Y3_{\, y3x3} - f_{x3}\0 \,\,= \,Y3_{\,y3x3} - f_{x3}\0,
\]
\[
f_{x3}\0 \,\,= \,X3_{\,x3x3} + X3_{\,y3x3}\0 \,\,= \,X3_{\,y3x3}\0. 
\]
These relations yield together 
\be\label{44}
f\,Y4_{\,x3}\0 \,\,= \,Y3_{\,y3x3} - X3_{\,y3x3}\0.
\ee
By the same reason as above the functions $X3$ and $Y3$ are affine in $y3$, 
and differentiating now the relations \eqref{42} and \eqref{43} sidewise 
with respect to $y3$ at 0, 
\[
f\,Y4_{\,y3}\0 \,\,= \,Y3_{\,x3y3} + Y3_{\,y3y3} - f_{y3}\0 \,\,= \,Y3_{\,x3y3} - f_{y3}\0,
\]
\[
f_{y3}\0 \,\,= \,X3_{\,x3y3} + X3_{\,y3y3}\0 \,\,= \,X3_{\,x3y3}\0, 
\]
which relations together imply 
\be\label{45}
f\,Y4_{\,y3}\0 \,\,= \,Y3_{\,x3y3} - X3_{\,x3y3}\0.
\ee
Yet the functions: $X3$ (available in level $X2$) and $Y3$ (available in level $Y2$) 
do not have second order terms $x3\,y3$, neither. So the quantities on the right hand 
sides in \eqref{44} and \eqref{45} vanish. Observation 1 is proved.\\ 

\n In view of Observation 1 the key relation \eqref{41} reduces to 
\be\label{46}
c\,Y4_{\,y4}\0 \,\,\,= \,{\t c}\,f\0.
\ee
{\bf Lemma.} $X3_{\,y3}\0 \,\,\,= \,0$.\\

\n Proof. The equation at level $X2$ reads in explicit terms 
\[
x2(\cdots) + \big(H + x2\,H_{\,x2}\big)\big(1 + x3\big) + x2\,y3\,H_{\,y2} \,= fG(1 + X3)
\]
Differentiating this equation sidewise with respect to $y3$ at 0 and remembering 
that neither $H$ nor $fG$ depends on $y3$ (cf. \eqref{28}), $0 = fG\,X3_{y3}\0$. 
Lemma is proved.\\

We are now in a position to replace $Y4$ in \eqref{46} by $Y3$. Namely, 
differentiating the relation \eqref{42} sidewise with respect to $y4$ at 0, 
\[
Y3_{\,y3}\0 \,\,= \,f_{\,y4} + f\,Y4_{\,y4}\0 \,\,= X3_{\,y3} + f\,Y4\0 \,\,= \,f\,Y4\0
\]
by Lemma. Whence the key relation \eqref{46} reduces to 
\be\label{47}
c\,Y3_{\,y3}\0 \,\,\,= \,{\t c}\,f^2\0.
\ee
{\bf Observation 2.} $f\0 \,\,= \,1$.\\

\n Proof. From \eqref{43} and Lemma there is 
\[
f\0 \,\,= \,X3_{\,x3} + X3_{\,y3}\0 \,\,= \, X3_{\,x3}\0\,.
\]
Let us write, for the reason of legibility only, 
$f\0 \,\,= \,f_0$, \,$fG\0 \,\,= \,(fG)_0$, \,$H\0 \,\,= \,H_0$. 
Then, in the Taylor expansion of the function $X3$, $X3 = f_0\,x3 + \cdots$, 
and also 
\be\label{48}
fG = (fG)_0 + {\rm h.o.t.}\,,\qquad\qquad H = H_0 + {\rm h.o.t.}
\ee
where in both the higher order terms above there is no coordinate $x3$ 
whatsoever (cf. also \eqref{28}). Let us write in explicit terms 
the equation in level $X2$
\be\label{49}
x2(\cdots) + \Big(H_0 + {\rm h.o.t.} + x2\,H_{\,x2}\Big)\big(1 + x3\big) + x2(\cdots) = 
\Big((fG)_0 + {\rm h.o.t.}\Big)\big(1 + f_0\,x3 + \cdots\big).
\ee
Now, by comparing the zero order terms in \eqref{49} 
\be\label{50}
H_0 = (fG)_0\,,
\ee
while comparing the $(x3)^1$ terms in \eqref{49} yields 
\be\label{51}
H_0 = (fG)_0\cdot f_0\,.
\ee
The relations \eqref{50} and \eqref{51} together yield $f_0 = 1$. 
Observation 2 is proved.\\

\n The relation \eqref{42} enhanced by Observation 2 clearly implies 
\be\label{52}
Y3_{\,x3} + Y3_{\,y3}\0 \,\,= \,1. 
\ee
{\bf Remark.} At this moment what only remains to show is that the first summand 
$Y3_{\,x3}\0$ in \eqref{52} vanishes. This follows from a rather compact series 
of inferences.\\

\n Firstly, level $Y1$ written in explicit terms 
\[
fG\cdot Y2 = x2\,(\ast) + Y1_{\,x1} + y2\,(\ast)
\]
implies that 
\be\label{53}
0 = fG\cdot Y2\0 \,\,= \,Y1_{\,x1}\0
\ee
and 
\be\label{54}
fG\cdot Y2_{\,x1}\0 \,\,= \,Y1_{\,x1x1}\0.
\ee
Secondly, level $Y0$ written in explicit terms 
\[
fGH\cdot Y1 = Y0_t + x1\,Y0_{\,x0} + y1\,Y0_{\,y0}
\]
implies 
\[
0 = \big(fGH\cdot Y1\big)_{\,x1x1}\0 \,\,= \,0 + 2\big(fGH\big)_{x1}\cdot Y1_{\,x1} + fGH\cdot Y1_{\,x1x1}\0 
\,\,= \,fGH\cdot Y1_{\,x1x1}\0
\]
by \eqref{53}. Hence 
\be\label{55}
Y2_{\,x1}\0 \,\,= \,0 
\ee
by \eqref{54}. Thirdly, level $Y2$ written in explicit terms 
\[
fG\cdot Y3 = x2\,(\ast) + Y2_{\,x1} + y2\,(\ast) + (1 + x3)Y2_{\,x2} + y3\,Y2_{\,y2}
\]
implies that 
\be\label{56}
0 = Y2_{\,x1} + Y2_{\,x2}\0
\ee
and also that 
\be\label{57}
fG\cdot Y3_{\,x3}\0 \,\,= \,Y2_{\,x2}\0.
\ee
Now \eqref{55} and \eqref{56} yield $Y2_{\,x2}\0 \,\,= \,0$. This latter 
equality coupled with \eqref{57} implies $Y3_{\,x3}\0 \,\,= \, 0$.\\
As it has been already noticed, this information ends the entire proof, 
implying by \eqref{52} that $Y3_{\,y3}\0 \,\,= \,1$. Hence, all in all, 
giving $c\cdot 1 = {\t c}\cdot 1$ (compare \eqref{47} and Observation 2).\\

Every two {\it different\,} values $c$ and $\t c$ of the parameter 
in the (pseudo) normal form \eqref{121} are non-equivalent. 

\section{Appendix}
We are going to demonstrate the techniques, standing behind the last statement 
in Theorem \ref{main}, on the singularity class 1.2.1.3.1. There are three pairwise 
non-equivalent types of EKR visualisations of the distribution germs sitting in 
that class, because its `father' class 1.2.1.3 is built up by three orbits 
(cf. section 4.2 in \cite{MP1}).\\

\n The members of the prolongation of the thickest orbit (of codimension 3, an open 
and dense part of 1.2.1.3) falling into the class 1.2.1.3.1 are visualised as the 
germs at $0 \in \R^{13}$ as follows 
\begin{align}\label{123one}
dx0 - x1dt         &= 0   &  dy0 - y1dt       &= 0\nonumber\\
dt   - x2\,dx1     &= 0   &  dy1 - y2\,dx1    &= 0\nonumber\\
dx2 - (1+x3)\,dx1  &= 0   &  dy2 - y3\,dx1    &= 0\nonumber\\
dx1 - x4\,dy3      &= 0   &  dx3 - y4\,dy3    &= 0\\
dx4 - (b+x5)\,dy3  &= 0   &  dy4 - (c+y5)\,dy3&= 0\nonumber\,,
\end{align}
where $b$ and $c$ are free real parameters. The cases of \,$b = 0$ or $c = 0$ are 
easy to handle -- to normalize by rescalings, respectively, the parameter $c$ to 1 
or $b$ to 1, unless $c$ or $b$ is 0. The (general) case $b \ne 0 \ne c$ is more 
interesting. The following rescaling of the EKR variables in \eqref{123one} 
(i.\,e., passing to new bar variables) simultaneously normalizes $b$ to 1 
and $c$ to ${\rm sgn}(c)$: 
\[
t = |c|^{-2}b^2\,\b{t}\,,\quad x0 = |c|^{-3}b^3\,\b{x0}\,,\quad y0 = |c|^{-9/2}b^4\,\b{y0}\,,
\quad x1 = |c|^{-1}b\,\,\b{x1}\,,
\]
\[
y1 = |c|^{-5/2}\,b^2\,\b{y1}\,,\quad x2 = |c|^{-1}b\,\,\b{x2}\,,\quad y2 = |c|^{-3/2}\,b\,\,\b{y2}\,,
\quad x3 = \b{x3}\,,
\]
\[
y3 = |c|^{-1/2}\,\b{y3}\,,\quad x4 = |c|^{-1/2}\,b\,\,\b{x4}\,,\quad y4 = |c|^{1/2}\,\b{y4}\,,
\quad x5 = b\,\b{x5}\,,\quad y5 = |c|\,\b{y5}\,.
\]
\vskip2mm
\n{\it Attention.} It is likely that the normal forms \eqref{123one} 
with $b = c = 1$ and with $b = 1$, $c = -1$ are not equivalent.\\

The members of the prolongations in 1.2.1.3.1 of the second orbit in 1.2.1.3 
(of codimension 4) are visualised in \eqref{121two} here below. We discuss now 
these pseudo-normal forms from the point of view of possible reductions of 
the parameters $b$ and $c$. 
\begin{align}\label{121two}
dx0 - x1dt         &= 0   &  dy0 - y1dt       &= 0\nonumber\\
dt   - x2\,dx1     &= 0   &  dy1 - y2\,dx1    &= 0\nonumber\\
dx2 - x3\,dx1      &= 0   &  dy2 - (1+y3)\,dx1&= 0\nonumber\\
dx1 - x4\,dy3      &= 0   &  dx3 - y4\,dy3    &= 0\\
dx4 - (b+x5)\,dy3  &= 0   &  dy4 - (c+y5)\,dy3&= 0\,.\nonumber
\end{align}
These pseudo-normal forms are particularly interesting because the computation 
of their i.s.'s by means of \cite{W} aroused initially a controversy. We reproduce 
here the value at $0 \in \R^{13}$ of the infinitesimal symmetry being the prolongation 
in the chart {\1}.{\2}.{\1}.{\3}.{\1} of an arbitrary smooth vector field $A\,\p_t + 
B\,\p_{x0} + C\,\p_{y0}$ in the base $\R^3(t,x0,y0)$. ($A,\,B,\,C$ are smooth functions 
of $t,\,x0,\,y0$, cf. Section 6 in \cite{MP2}.) The outcome of the computation depends 
only on the first jets of $A,\,B,\,C$ at $0 \in \R^3$: 
\begin{align}\label{ivec}
A\,\p_t + &B\,\p_{x0} + C\,\p_{y0} + B_t\p_{x1} + C_t\p_{y1} + C_{x0}\p_{y2} 
+ \big(A_t - 2B_{x0} + C_{y0}\big)\p_{y3}\nonumber\\
+ &\,b\,\big(-3A_t + 5B_{x0} - 2\,C_{y0}\big)\p_{x5} + c\,\big(A_t + 2B_{x0} 
- 2\,C_{y0}\big)\p_{y5}\0\,.
\end{align}
The idea of using the i.s.'s in the classification problem consists in assuming 
the vanishing at $0 \in \R^{13}$ of all the components of \eqref{ivec} save 
possibly these at the newly incoming directions $\p_{x5}$ and $\p_{y5}$. 
Then, in the case of the $\p_{x5}$- and/or $\p_{y5}$-components turning out 
to be automatically zero, the variables $x5$ and/or $y5$ cannot be moved from 
their initial values $b$ and/or $c$ by means of the automorphisms of $\D^5$ 
embeddable in flows, keeping all the older variables frozen. That would be a hint 
that the parameter $b$ and/or $c$ {\it could\,} be a modulus of the local classification 
-- to be verified separately.\\
However, the expression \eqref{ivec} (got after rectifying a preliminary flawed output) 
is {\it not\,} of that kind. Assuming the vanishing of the components of \eqref{ivec} 
save the last two, the shorter expression of that vector is 
\be\label{ivect}
\big(-B_{x0} + C_{y0}\big)\,b\,\p_{x5} + \big(4B_{x0} - 3\,C_{y0}\big)\,c\,\p_{y5}\0\,.
\ee
For $b \ne 0 \ne c$ one can vary the functions $B$ and $C$ so as to get in \eqref{ivect} 
non-zero coefficients at $b\,\p_{x5}$ and $c\,\p_{y5}$ one independently of the other. 
Hence a modulus of the local classification is excluded in this case, although 
the precise number of orbits is not ascertained at this moment.\\

\n In parallel, elementary rescalings show that the germs \eqref{121two} belong 
to a finite number of orbits. In fact (as previously), the cases $b = 0$ and $c = 0$ 
are particularly easy to handle -- to normalize, respectively, the parameter $c$ 
to 1 and $b$ to 1, unless $c$ or $b$ is 0.\\
In the general case $b \ne 0 \ne c$ the following rescaling of variables 
normalizes both $b$ and $c$ to 1: 
\[
t \,= \,b^2c\,\,\b{t}\,,\quad x0 \,= \,b^3c\,\,\b{x0}\,,\quad y0 \,= \,b^4c\,\,\b{y0}\,,
\quad x1 \,= \,b\,\,\b{x1}\,,
\]
\[
y1 \,= \,b^2\,\b{y1}\,,\quad x2 \,= \,b\,c\,\,\b{x2}\,,\quad y2 \,= \,b\,\,\b{y2}\,,
\quad x3 \,= \,c\,\,\b{x3}\,,
\]
\[
y3 \,= \,\b{y3}\,,\quad x4 \,= \,b\,\,\b{x4}\,,\quad y4 \,= \,c\,\,\b{y4}\,,
\quad x5 \,= \,b\,\,\b{x5}\,,\quad y5 \,= \,c\,\,\b{y5}\,.
\]
\vskip2mm
Endly, the members of the prolongations in 1.2.1.3.1 of the third orbit 
in 1.2.1.3 (of codimension 5) are visualised in \eqref{121three} here below. 
\begin{align}\label{121three}
dx0 - x1dt         &= 0   &  dy0 - y1dt       &= 0\nonumber\\
dt   - x2\,dx1     &= 0   &  dy1 - y2\,dx1    &= 0\nonumber\\
dx2 - x3\,dx1      &= 0   &  dy2 - y3\,dx1    &= 0\nonumber\\
dx1 - x4\,dy3      &= 0   &  dx3 - y4\,dy3    &= 0\\
dx4 - (b+x5)\,dy3  &= 0   &  dy4 - (c+y5)\,dy3&= 0\,.\nonumber
\end{align}
The discussion of the status of these pseudo-normal forms from the point of 
view of possible reductions of the parameters $b$ and $c$ is now extremely short. 
The subcases $b = 0$ or $c = 0$ are quickly normalizable, meaning the normalization 
to 1 of the other parameter, by simple rescalings. The leading subcase $b \ne 0 \ne c$ 
is normalizable (both $b$ and $c$ reducible to 1) by means of {\it the same} rescaling 
of variables as for \eqref{121two} above -- because that rescaling leaves 
the variable $y3$ untouched.\\

\n The last statement in Theorem \ref{main} is now partly justified 
on just one chosen example of the singularity class 1.2.1.3.1\,.\\

\n Acknowledgments. The author is grateful to Andrzej Weber for his help 
in carrying through the implementation of the earlier known recursive formulae 
and then in the effective computations of the infinitesimal symmetries of 
special 2-flags. 


\end{document}